\documentclass[12pt]{article}
\usepackage[utf8]{inputenc}
\usepackage[T2A]{fontenc}
\usepackage[russian,english]{babel}
\usepackage{amsmath,amsfonts,amssymb,euscript,graphicx,wrapfig,multirow}
\usepackage{dsfont}
\usepackage{amsthm}
\inputencoding{utf8}
\usepackage{hyperref}
\usepackage{cite}

\makeatletter
\renewcommand{\@biblabel}[1]{#1.}
\makeatother

\theoremstyle{theorem}

\theoremstyle{dfn}
\newtheorem{dfn}{Definition}
\theoremstyle{remark}

\begin{document}
%\renewcommand\refname{\centering References}

% Переключаем язык на английский.
% Очень полезно как в плане типографики (в том числе расстановки переносов),
% так и в плане того, что не надо переименовывать "Рисунок" в "Figure"
\selectlanguage{english}

\title{
	On particular examples of planar integral point sets and their classification
	\footnote{
		This work was carried out at Voronezh State University and supported by the Russian Science
		Foundation grant 19-11-00197.
	}
}

%% Прекрасно понимаю, что следующая команда - дичь и вордописчество, но время не ждёт, время жмёт
\author{
	Avdeev N.N.
	\footnote{nickkolok@mail.ru, avdeev@math.vsu.ru}
	, Momot E.A., Zvolinskiy A.E.
	\\
	\\
	\emph{Voronezh State University}
}

\maketitle

\paragraph{Abstract.}
	A planar integral point set is a set of non-collinear points in plane
	such that for any pair of the  points the Euclidean distance
	between them is integral.
	We discuss the classification of planar integral point sets
	and provide examples of sets that are not covered by the existent classification.

\paragraph{Keywords.}
	integral point sets,
	planar integral point sets,
	classification of planar integral point sets,
	diameter of an integral point set

\section{Introduction}

\begin{dfn}\label{dfn1}
	A planar integral point set (PIPS) is a set $\mathcal{P}$
	of non-collinear points in plane $\mathbb{R}^{2}$ such that
	for any pair of points $P_{1}, P_{2} \in \mathcal{P}$
	the Euclidean distance $|P_{1}P_{2}|$
	between points $P_{1}$ and $P_{2}$ is integral.
\end{dfn}

How do we describe an integral point set?
For example, by the number of points in it, which is always finite~\cite{anning1945integral,erdos1945integral}
and is said to be the \emph{cardinality} of the IPS.
Furthermore, we can naturally % TODO: naturally/essentially? cf. arxiv/1
define the diameter of a finite point set.

\begin{dfn}
	The diameter of the integral point set $\mathcal{P}$ is defined as
	\begin{equation}
		\operatorname{diam(\mathcal{P})} = \underset{P_{1}, P_{2} \in
		\mathcal{P}}{\max} |P_{1}P_{2}|
		.
	\end{equation}
\end{dfn}

Every planar IPS aso has a characteristic~\cite{kemnitz1988punktmengen,kurz2005characteristic}.

\begin{dfn}
	The characteristic of a planar integral point set $M$ is the least positive integer $q$
	such that the area of any triangle with vertices $M_1, M_2, M_3 \in M$
	is commensurable with $\sqrt{q}$.
\end{dfn}
Characteristic of an IPS does not change if the set is moved, dilated ot flipped;
moreover, even addition or deletion of a point does not change the characteristic of IPS.

While the minimal possible diameter for planar integral point sets of given cardinality was being computed,
it was noticed~\cite{kurz2008minimum} that such diameter is attained at sets with many points on a straight line;
for some estimations on this tendency, we also refer the reader to~\cite{solymosi2003note}.
Thus, the following classification was introduced:
\begin{dfn}
	A planar integral point set $M$ is said to be in \emph{semi-general position}
	if no three points of $M$ are located in a straight line.
\end{dfn}

The most dominating examples of PIPS in semi-general position are circular sets.
\begin{dfn}
	A planar integral point sets that is situated on a circle is said to be a \textit{circular}
	point set.
\end{dfn}

So, the following constraint appeared.
\begin{dfn}
	A planar integral point set $M$ is said to be in \emph{general position}
	if no three points of $M$ are located in a straight line
	and no four points of $M$ are located in a circle.
\end{dfn}

Planar integral point sets in general position are very difficult to find;
the first known examples were presented in~\cite{kreisel2008heptagon}.
As for now, there is no known example of PIPS of cardinality 8 in general position.

The main purpose of this work is to provide examples of planar integral point sets
that may give the clue for development of further classification.

For convenience, we use the notation \cite{our-ped-2018,our-pmm-2018,our-vmmsh-2018}:
$\sqrt{p}/q * \{ (x_1,y_1), ...,$ $ (x_n, y_n)  \}$,
which means that each abscissa is multiplied by $1/q$
and each ordinate is multiplied by $\sqrt{p}/q$,  i.e.
\begin{equation}
	\label{eq:char_lattice}
	\sqrt{p}/q * \{ (x_1,y_1), ..., (x_n, y_n)  \}
	=
	\left\{ \left(\frac{x_1}{q},\frac{y_1\sqrt{p}}{q}\right), ..., \left(\frac{x_n}{q},   \frac{y_n\sqrt{p}}{q}\right)  \right\}
	.
\end{equation}
Here $q$ is the characteristic of the PIPS;
every PIPS can be represented in such way~\cite{our-vmmsh-2018-translit}[Theorem 4].

Note that all examples that are discussed below are
located on a union of at most three straight lines.
%In the most pictures these lines are shown as solid lines.
For classification of planar IPS that are located on a union of two straight lines,
we refer the reader to~\cite{avdeev2019particular}.

There are some examples of planar IPS that are not contained in the union of any three straight lines:
for examples, these are heptagons presented in~\cite{kreisel2008heptagon} and 7-clusters from~\cite{kurz2013constructing}.
However, we have to keep in mind that the circular inversion under certain conditions
translates an integral point set into an integral point set
(although sometimes additional dilation is necessary).
On the other hand, the circular inversion may translate a straight line into a circle and vice versa.
Thus, we can consider \emph{generalized circles}, that are circles or straight lines;
obviously, in that point of view all the examples from~\cite{kreisel2008heptagon} and~\cite{kurz2013constructing}
are located on a union of three generalized circles,
because each seven points are located on a union of a circle and two straight lines.

\section{Rails sets}

\begin{dfn}
	A planar integral point sets of $n$ points with $n-1$ points on a straight line is called
	a \textit{facher} set.
\end{dfn}
Facher sets are predominant examples of planar integral point sets.
In~\cite{antonov2008maximal}, facher sets of characteristic 1 are called \textit{semi-crabs}.
%For $9 \leq n \leq 122$, the diameter $d(2,n)$ is reached on a facher point set~\cite{kurz2008minimum}.

\begin{dfn}[\cite{avdeev2019particular}]
	A non-facher planar integral point sets situated in two parallel straight lines
	is called a \textit{rails} set.
\end{dfn}

Among the rails sets, sets with 2 points on one line and all the other on another line dominate.

Two IPSs below have been obtained by dilating~\cite[Fig. 34]{avdeev2019particular} by $23$ and $29$  resp.;
third one has been constructed by dilation and merge.

\begin{figure}[h!]
\center{\includegraphics[width=1\linewidth]{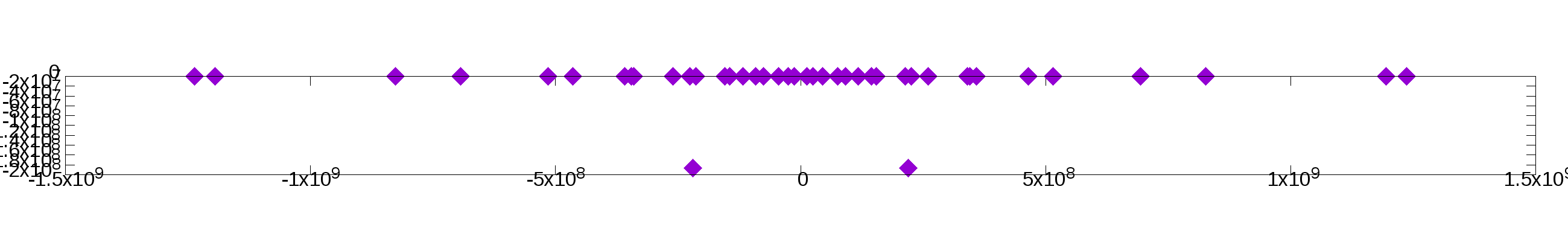}}
\parbox{1\linewidth}{\caption{IPS of cardinality 42 and diameter 2473117504}
\label{42_symm.png}}
\end{figure}

\begin{figure}[h!]
\center{\includegraphics[width=1\linewidth]{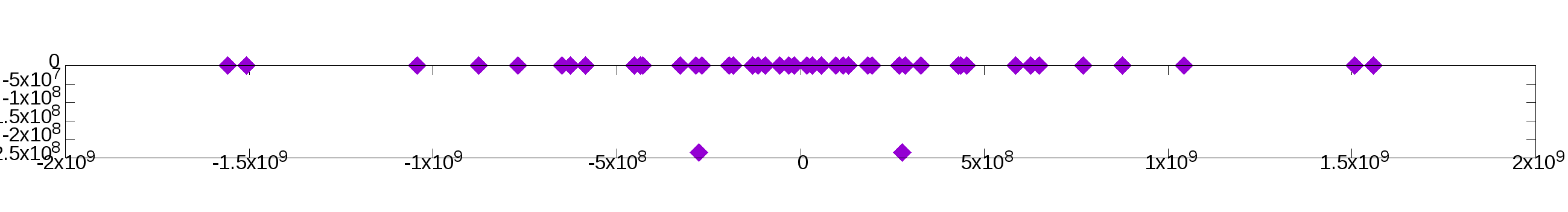}}
\parbox{1\linewidth}{\caption{IPS of cardinality 46 and diameter 3118278592}
\label{46_symm.png}}
\end{figure}

\begin{figure}[h!]
\center{\includegraphics[width=1\linewidth]{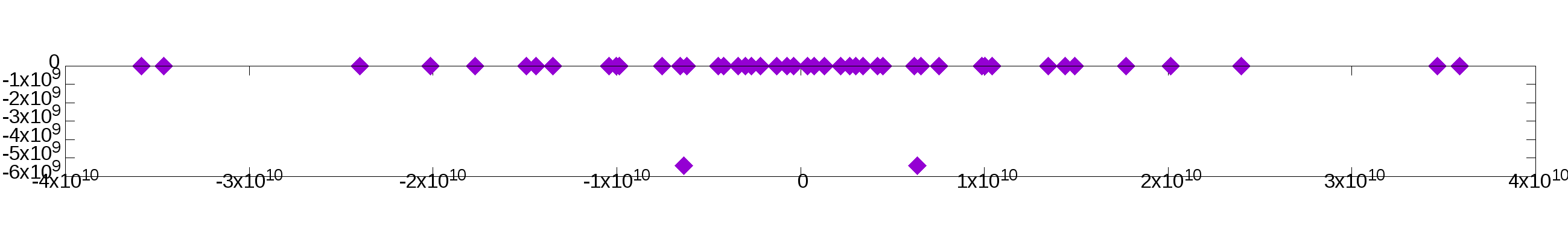}}
\parbox{1\linewidth}{\caption{IPS of cardinality 48 and diameter 71720407616}
\label{48_symm.png}}
\end{figure}

\begin{itemize}
\setlength{\itemsep}{-1mm}

\item
Figure~\ref{42_symm.png}:
\begin{multline}
	\mathcal{P}_{42}=\sqrt{154}/{1} * \{
		(\pm219513840; -15069600);
		(\pm 345596160; 0);
		(\pm260201760; 0);
		\\
		(\pm225792840; 0);
		(\pm213234840; 0);
		(\pm153961080; 0);
		(\pm144668160; 0);
		(\pm25116000; 0);
		\\
		(\pm694026840; 0);
		(\pm514710560; 0);
		(\pm359116940; 0);
		(\pm13423904; 0);
		(\pm75682880; 0);
		\\
		(\pm464143680; 0);
		(\pm827069880; 0);
		(\pm92144325; 0);
		(\pm1195180740; 0);
		(\pm1236558752; 0);
		\\
		(\pm44590560; 0);
		(\pm339925740; 0);
		(\pm117312468; 0)
	\}
\end{multline}

\item
Figure~\ref{46_symm.png}:
\begin{multline}
	\mathcal{P}_{46}=
	\sqrt{154}/{1} * \{
		(\pm276778320; -19000800);
		(\pm435751680; 0);
		(\pm328080480; 0);
		\\
		(\pm268861320; 0);
		(\pm194124840; 0);
		(\pm182407680; 0);
		(\pm1559139296; 0);
		(\pm284695320; 0);
		\\
		(\pm1506967020; 0);
		(\pm1042827240; 0);
		(\pm875077320; 0);
		(\pm648982880; 0);
		(\pm585224640; 0);
		\\
		(\pm452799620; 0);
		(\pm95426240; 0);
		(\pm16925792; 0);
		(\pm116181975; 0);
		(\pm428602020; 0);
		\\
		(\pm56222880; 0);
		(\pm769560480; 0);
		(\pm626458560; 0);
		(\pm31668000; 0);
		(\pm130761918; 0)
	\}
\end{multline}

\item
Figure~\ref{48_symm.png}:
\begin{multline}
	\mathcal{P}_{48}=
	\sqrt{154}/{1} * \{
		( \pm6365901360 ; -437018400);
		( \pm10022288640; 0);
		( \pm23985026520 ; 0);
		\\
		( \pm389293216 ; 0);
		( \pm6183810360 ; 0);
		( \pm4464871320 ; 0);
		( \pm4195376640 ; 0);
		( \pm728364000 ; 0);
		\\
		( \pm35860203808 ; 0);
		( \pm34660241460 ; 0);
		( \pm7545851040 ; 0);
		( \pm20126778360 ; 0);
		\\
		( \pm14926606240 ; 0);
		( \pm13460166720 ; 0);
		( \pm10414391260 ; 0);
		( \pm2194803520 ; 0);
		\\
		( \pm2672185425 ; 0);
		( \pm9857846460 ; 0);
		( \pm1293126240 ; 0);
		( \pm17699891040 ; 0);
		\\
		( \pm14408546880 ; 0);
		( \pm3007524114 ; 0);
		( \pm6547992360 ; 0);
		( \pm3402061572 ; 0)
	\}
\end{multline}

\end{itemize}

Taking the examples into consideration,
we can conjecture that there is an infinite point set with rational distances
that contains $\mathcal{P}_{48}$.
(However, it is known~\cite{solymosi2010question} that
if a point set $S$ with rational distances has infinitely many points on a line or on a circle,
then all but 4 resp. 3 points of $S$ are on the line or on the circle.)

\section{Example of sets with many common points that cannot be merged}

\begin{figure}[h!]
	\begin{minipage}[h]{0.32\linewidth}
		\begin{center}
			\includegraphics[width=1\linewidth]{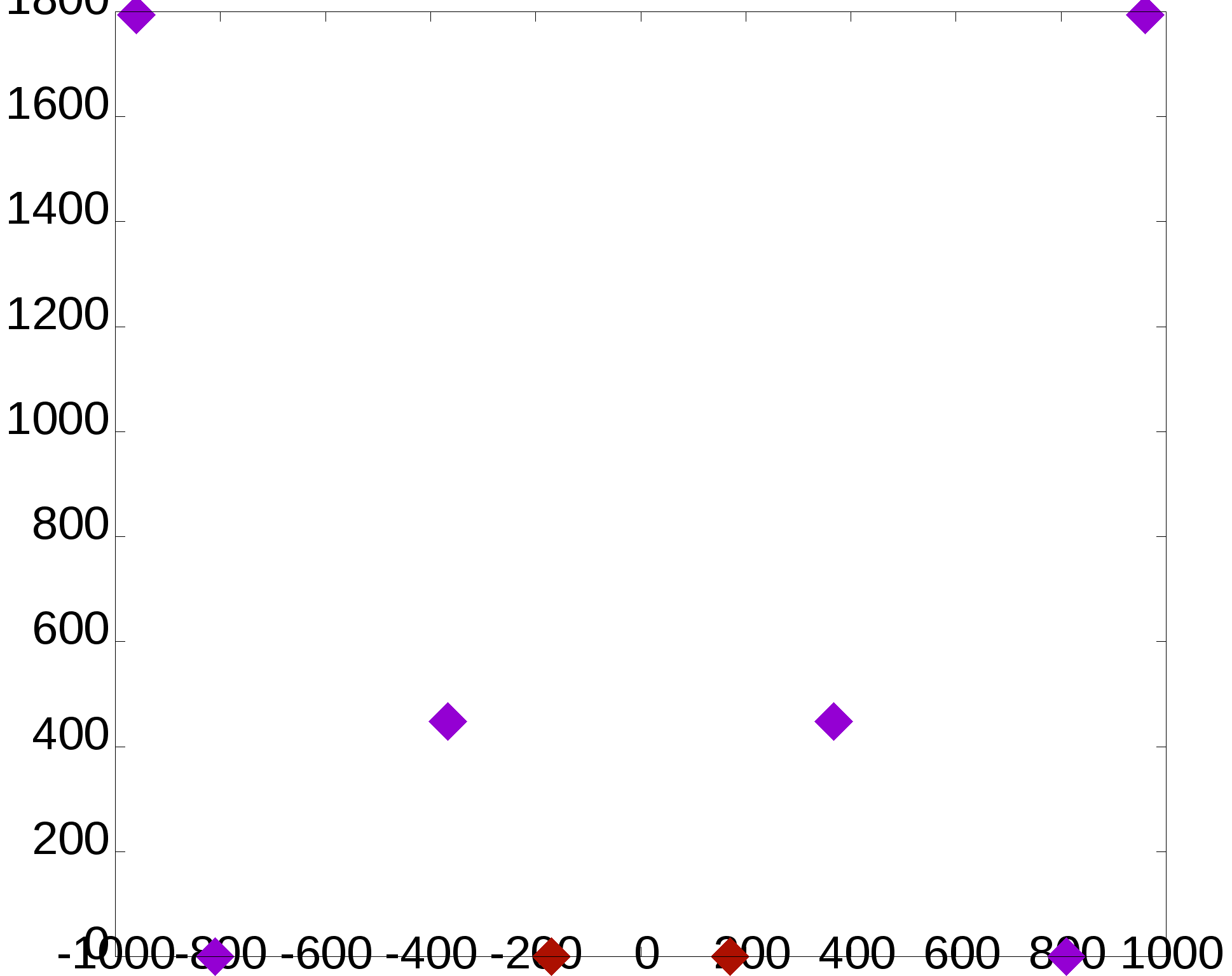}\\ a)
		\end{center}
	\end{minipage}
	\hfill
	\begin{minipage}[h]{0.32\linewidth}
		\begin{center}
			\includegraphics[width=1\linewidth]{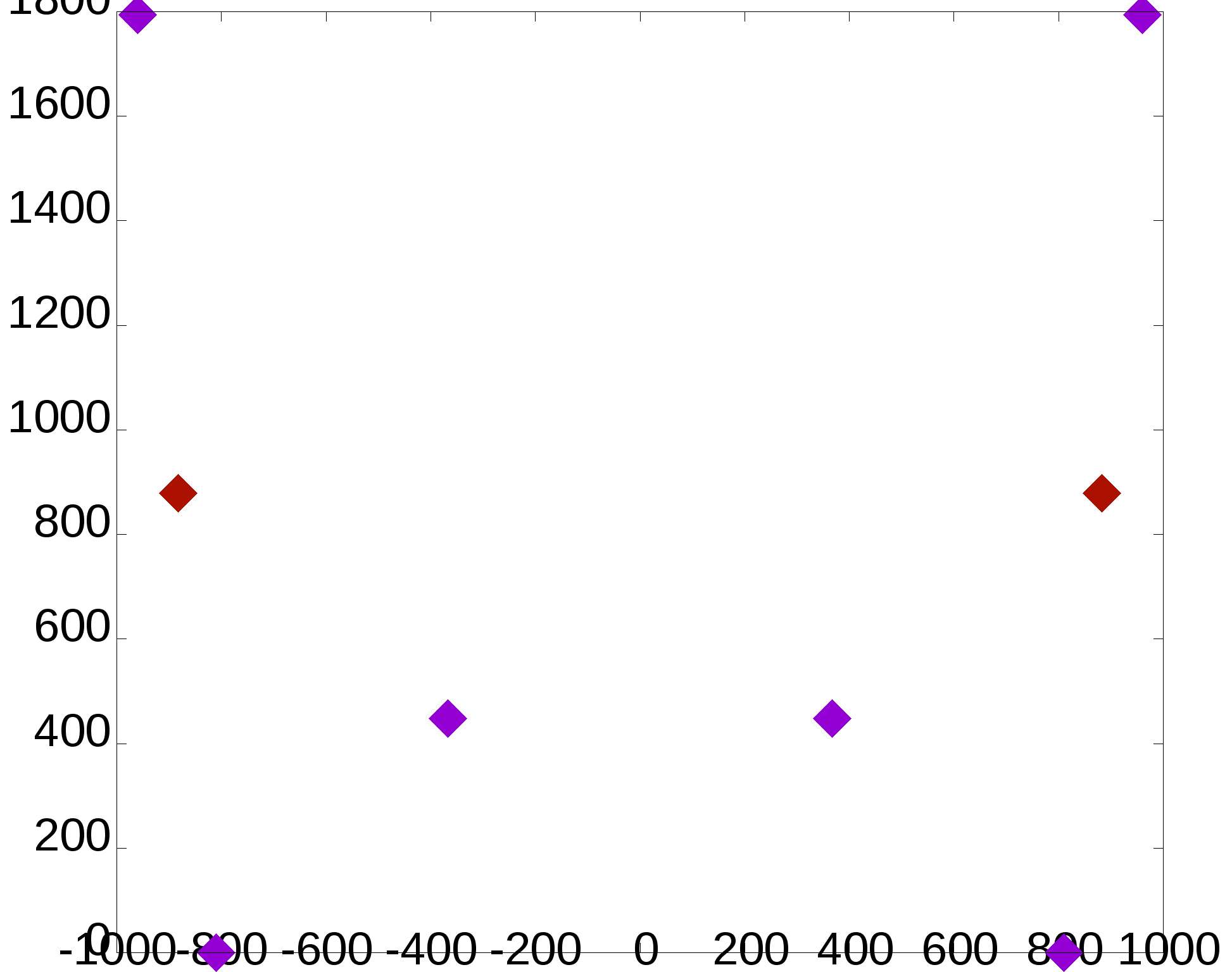}\\ b)
		\end{center}
	\end{minipage}
	\begin{minipage}[h]{0.32\linewidth}
		\begin{center}
			\includegraphics[width=1\linewidth]{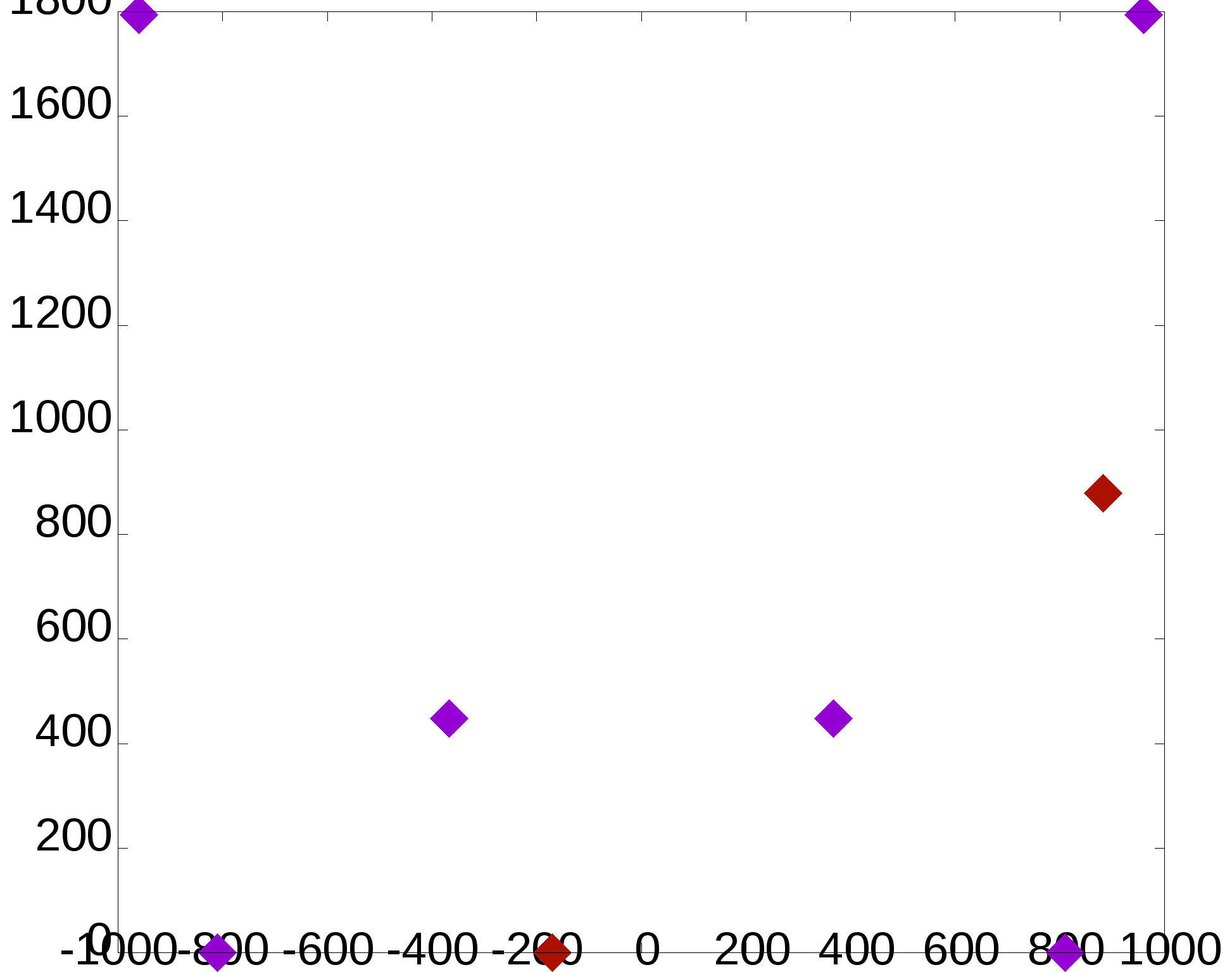}\\ c)
		\end{center}
	\end{minipage}
	\hfill
	\caption{IPSs with cardinality 8 and diameter 2520 with many common points}
	\label{8_with_many_common}
\end{figure}

Figure~\ref{8_with_many_common} shows an example of three PIPS of cardinality 8,
each pair of that shares 6 or 7 points but cannot be combined into another PIPS.

\begin{itemize}
\item
$\mathcal{P}=\sqrt{143}/2*\{
( \pm1620 ; 0);
( \pm1920 ; 300);
( \pm735 ; 75);
( \pm340 ; 0);
\}$

\item
$\mathcal{P}=\sqrt{143}/2*\{
( \pm1620 ; 0);
( \pm1920 ; 300);
( \pm735 ; 75);
( \pm1767 ; 147);
\}$

\item
$\mathcal{P}=
\sqrt{143}/2*\{
( \pm1620 ; 0);
( \pm1920 ; 300);
( \pm735 ; 75);
( -340 ; 0);
( 1767 ; 147);
\}$

\end{itemize}

The distance between the non-adoptable points is
\begin{equation}
	\sqrt{\left(\frac{1767}{2} - \frac{340}{2}\right)^2 + \left(\frac{147}{2}\right)^2\cdot143} = 2\sqrt{320401}
	.
\end{equation}
It's notable that 320401 is a prime number.

\section{Integral point sets with two axes of symmetry}

\begin{figure}[h!]
\center{\includegraphics[width=1\linewidth]{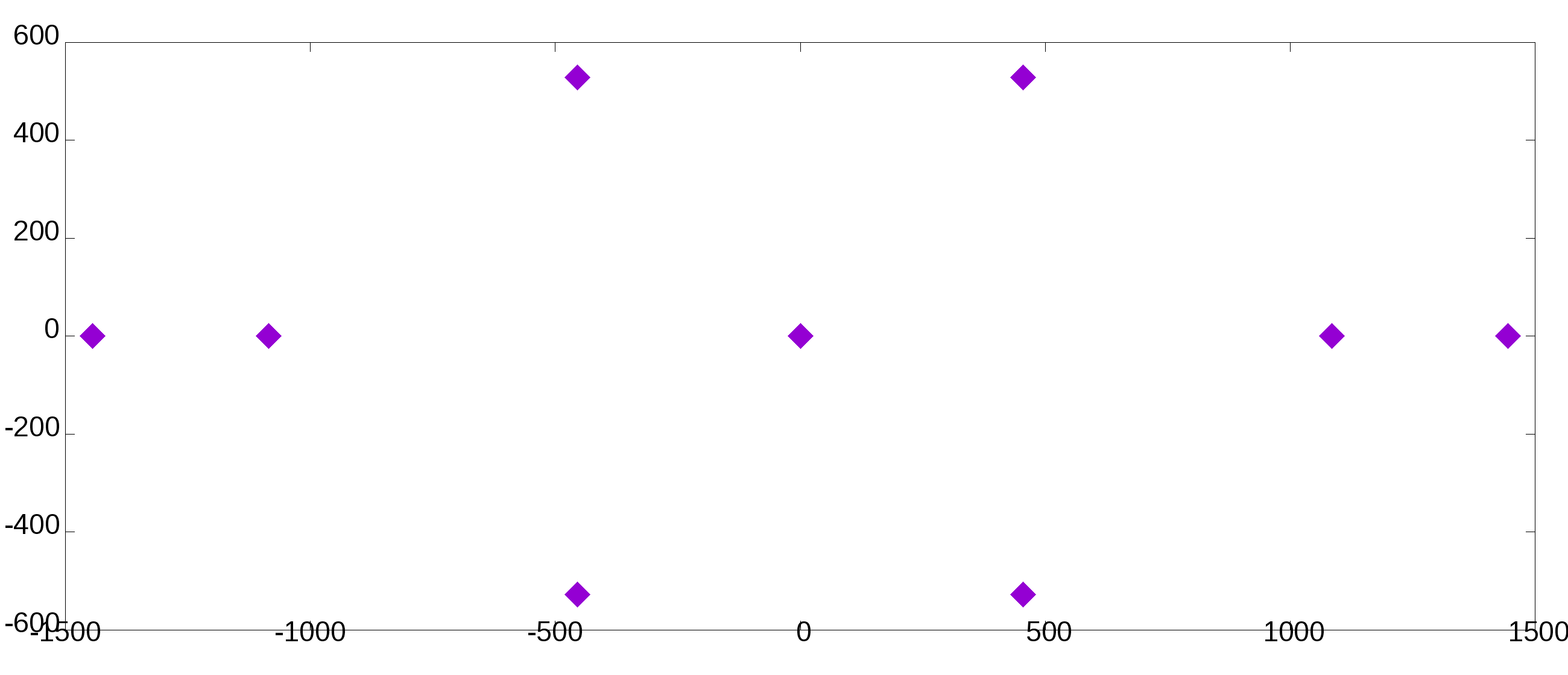}}
\parbox{1\linewidth}{\caption{IPS of cardinality 9 and diameter 2890}
\label{9_2890_1_02e6af}}
\end{figure}

\begin{figure}[h!]
\center{\includegraphics[width=1\linewidth]{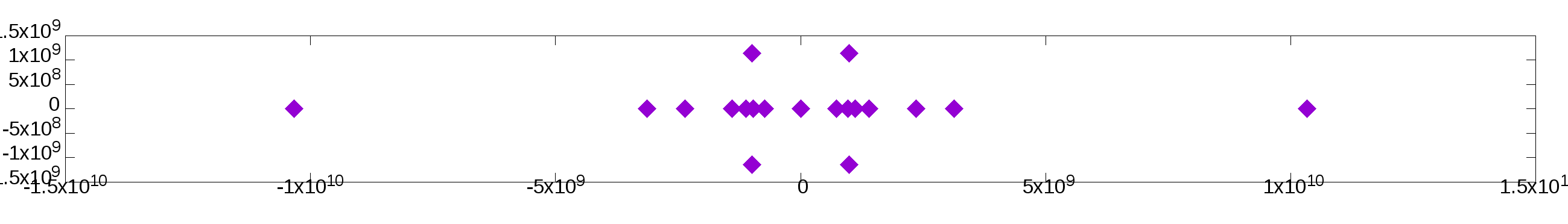}}
\parbox{1\linewidth}{\caption{IPS of cardinality 19 and diameter 20663808074}
\label{19_20663808074}}
\end{figure}

The set $\mathcal{P}_{19}$ shown on Figure~\ref{19_20663808074} was obtained from the set $\mathcal{P}_{9}$ shown on Figure~\ref{9_2890_1_02e6af}
by dilation and looking for points on $x$ axis.

\begin{equation}
	\mathcal{P}_9=\sqrt{1}/1*\{
	( 0 ; 0);
	( \pm1445 ; 0);
	( \pm1085 ; 0);
	( -455 ; \pm528);
	( 455 ; \pm528)
	\}
\end{equation}

\begin{multline*}
	\mathcal{P}_{19}=
	\sqrt{1}/1*\{
		( 0 ; 0);
		( -987843675 ; \pm1146332880);
		( 987843675 ; \pm1146332880);
		\\
		( \pm729918777 ; 0);
		( \pm972103809 ; 0);
		( \pm1113030324 ; 0);
		( \pm1400170149 ; 0);
		\\
		( \pm3137217825 ; 0);
		( \pm2355627225 ; 0);
		( \pm10331904037 ; 0);
	\}
\end{multline*}

\section{Arrow-like integral point sets with one axis of symmetry}

\begin{figure}[h!]
\center{\includegraphics[width=1\linewidth]{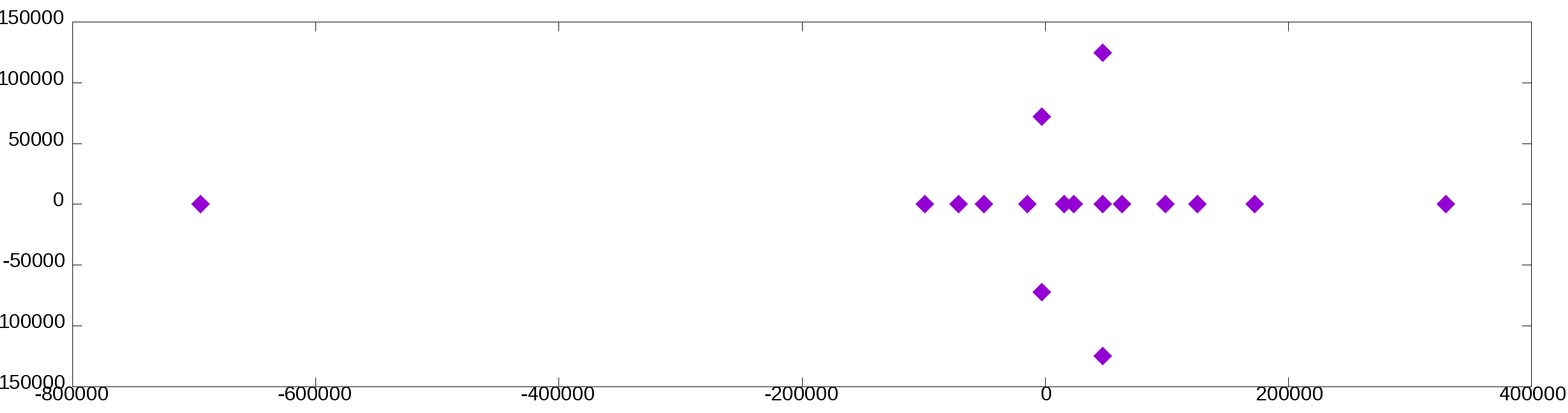}}
\parbox{1\linewidth}{\caption{IPS of cardinality 17 and diameter 1024296}
\label{17_1024296_1_639b}}
\end{figure}

On Figure~\ref{17_1024296_1_639b}, the following IPS is shown:
\begin{multline}
	\mathcal{P}_{17}=
	\sqrt{1}/1*\{
		( -2847 ; \pm72072);
		( 47073 ; \pm124488);
		( 47073 ; 0);
		( \pm98943 ; 0);
		\\
		( -694668 ; 0);
		( -71487 ; 0);
		( -50367 ; 0);
		( -14943 ; 0);
		( 15457 ; 0);
		\\
		( 23582 ; 0);
		( 63073 ; 0);
		( 125307 ; 0);
		( 172207 ; 0)
		( 329628 ; 0);
	\}
\end{multline}
Note that the axis of symmetry for $\mathcal{P}_{17}$ is the $x$ axis;
all such sets are of characteristic 1.
Moreover, note that the set contains three points with the same first coordinates.

\begin{figure}[h!]
	\center{\includegraphics[width=1\linewidth]{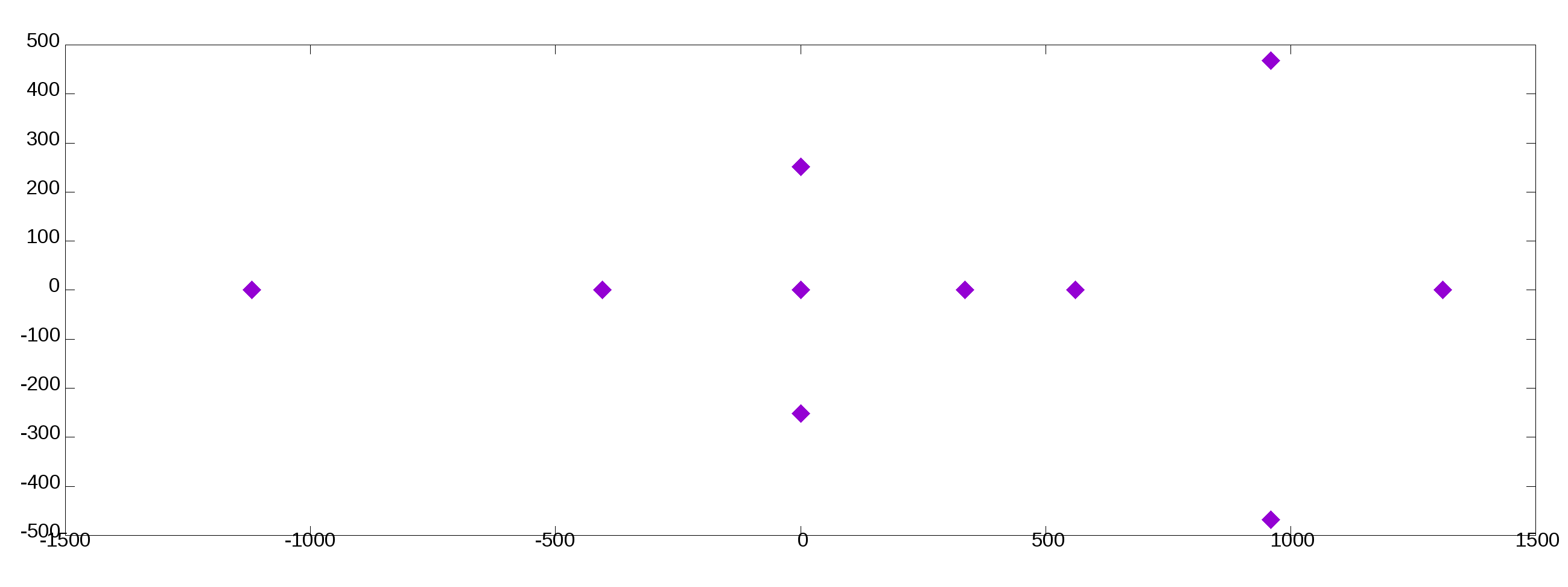}}
	\parbox{1\linewidth}{\caption{IPS of cardinality 10 and diameter 2431}
	\label{10_2431_1_0606}}
\end{figure}

\begin{figure}[h!]
	\center{\includegraphics[width=1\linewidth]{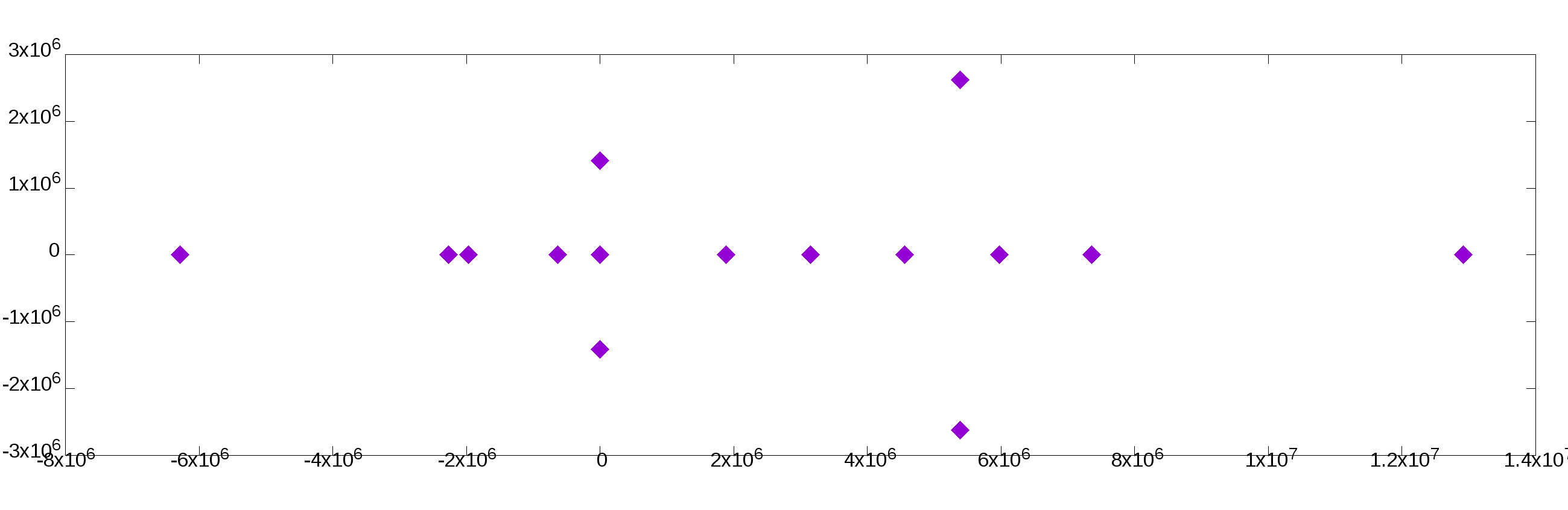}}
	\parbox{1\linewidth}{\caption{IPS of cardinality 15 and diameter 19203744}
	\label{15_19203744_1_80db}}
\end{figure}

On Figures~\ref{10_2431_1_0606} and~\ref{15_19203744_1_80db}, other examples of arrow-like IPSs are shown.
The one on Figure~\ref{15_19203744_1_80db} is obtained from the one on Figure~\ref{10_2431_1_0606}
by dilation and looking for points on $x$ axis.

Figure~\ref{10_2431_1_0606}:
\begin{multline*}
	\mathcal{P}_{10} =
	\sqrt{1}/1*\{
		( 0 ; 0);
		( 0 ; \pm252);
		( 960 ; \pm468);
		( -1120 ; 0);
		( -405 ; 0);
		( 336 ; 0);
		( 561 ; 0);
		( 1311 ; 0)
	\}
\end{multline*}

Figure~\ref{15_19203744_1_80db}:
\begin{multline*}
	\mathcal{P}_{15} =
	\sqrt{1}/1*\{
		( 0 ; 0);
		( 0 ; \pm1413720);
		( 5385600 ; \pm2625480);
		\\
		( -6283200 ; 0);
		( -2272050 ; 0);
		( -1971915 ; 0);
		( -635040 ; 0);
		( 1884960 ; 0);
		\\
		( 3147210 ; 0);
		( 4558176 ; 0);
		( 5976333 ; 0);
		( 7354710 ; 0);
		( 12920544 ; 0);
	\}
\end{multline*}

\section{Other examples}

\begin{figure}[h!]
\center{\includegraphics[width=1\linewidth]{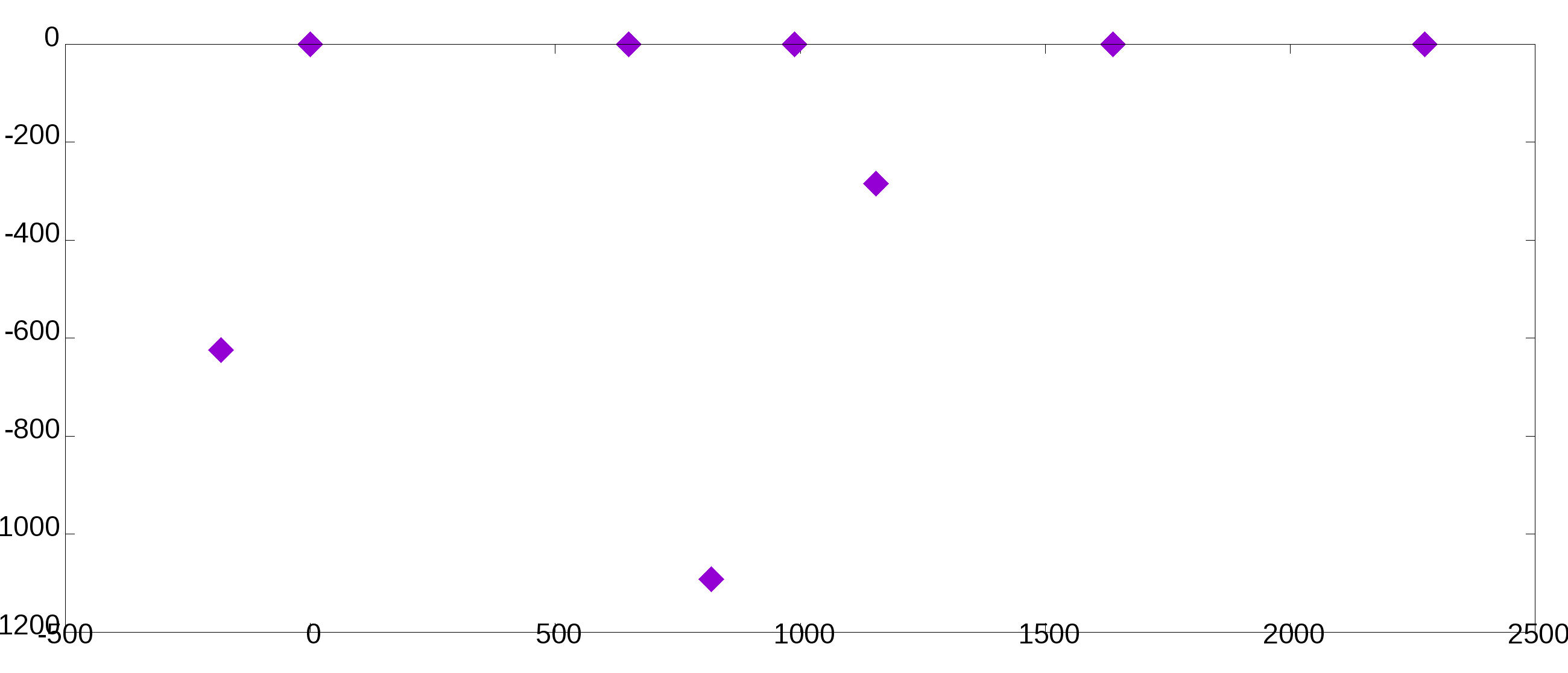}}
\parbox{1\linewidth}{\caption{IPS of cardinality 8 and diameter 2535}
\label{8_2535_1_d680}}
\end{figure}

Figure~\ref{8_2535_1_d680} displays an IPS with
no axis of symmetry;
although the set is of characteristic 1,
it cannot be extended by the reflection in the $x$ axis.
Moreover, we failed to extend it by dilation and looking for extra points on the $x$ axis.
%
%TODO: write about points which prevent the set from reflection
%
\begin{multline}
	\mathcal{P}_8=
	\sqrt{1}/13*
	\{
	( 0 ; 0);
	( 8450 ; 0);
	( 12844 ; 0);
	( 21294 ; 0);
	( 29575 ; 0);
	\\
	( -2366 ; -8112);
	( 10647 ; -14196);
	( 15022 ; -3696)
	\}
\end{multline}

\begin{figure}[h!]
	\center{\includegraphics[width=1\linewidth]{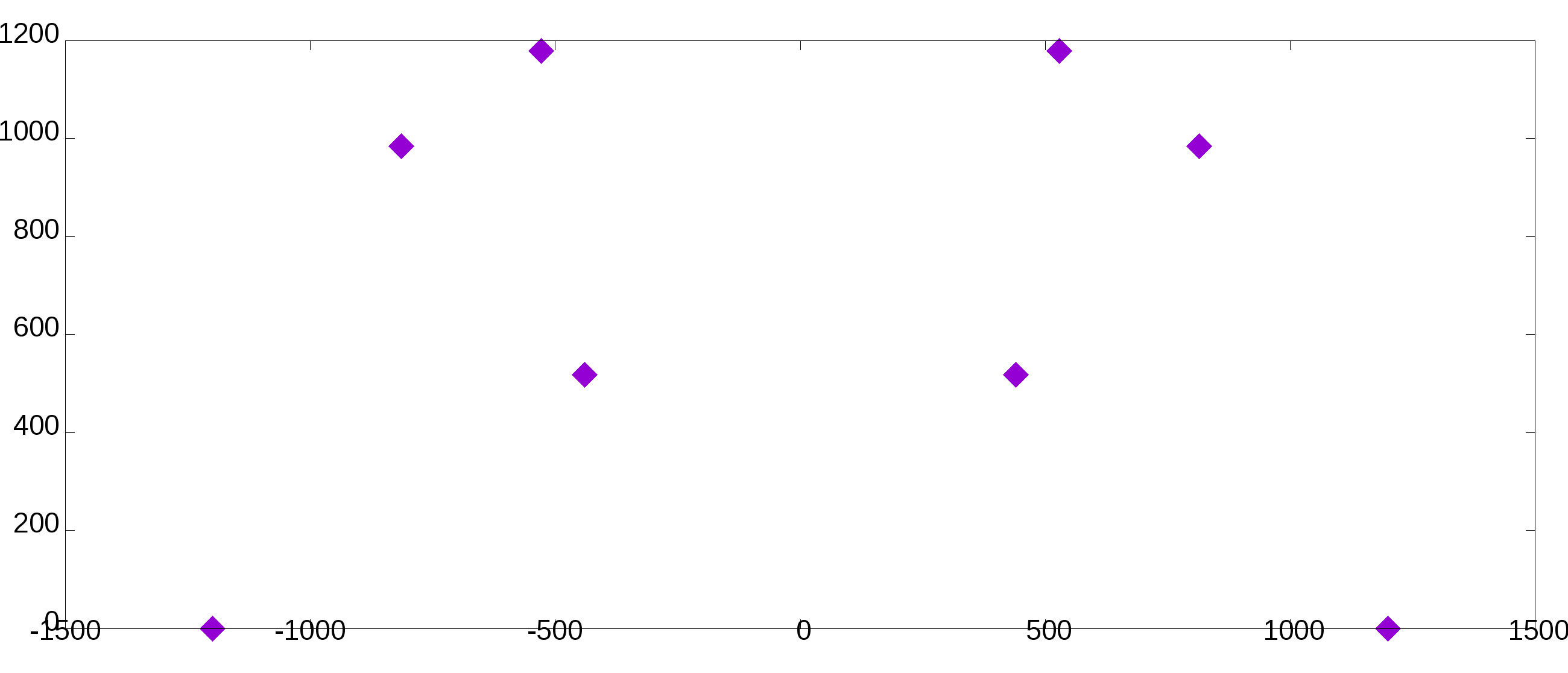}}
	\parbox{1\linewidth}{\caption{IPS of cardinality 8 and diameter 2400}
	\label{8_2400_42_56f3}}
\end{figure}

The set shown on Figure~\ref{8_2400_42_56f3} has an axis of symmetry, but it is $y$ axis, not $x$ axis.
Due to the fact that its characteristic is not 1,
the set cannot be rotated by $90^\circ$ but still be on lattice~\eqref{eq:char_lattice}:
\begin{equation}
	\mathcal{P}_{8y}=
	\sqrt{42}/1*\{( \pm1200 ; 0);
	( \pm529 ; 182);
	( \pm814 ; 152);
	( \pm440 ; 80)
	\}
\end{equation}

\section{Final remarks}
All the given planar integral point sets were obtained through a combination of computer search and intuition of the authors.

The source code can be obtained at https://gitlab.com/Nickkolok/ips-algo

\section{Acknowledgements}
The authors thank Dr. Prof. E.M. Semenov and Dr.~A.S.~Usachev for proofreading.

%\printbibliography
%\printbibitembibliography

The authors declare that they have no known competing financial interests or personal relationships that could have appeared to influence the work reported in this paper.

\end{document}